# Empirical processes indexed by estimated functions

## Aad W. van der Vaart[1] and Jon A. Wellner[2],*

*Vrije Universiteit Amsterdam and University of Washington*

**Abstract:** We consider the convergence of empirical processes indexed by functions that depend on an estimated parameter $\eta$ and give several alternative conditions under which the "estimated parameter" $\eta_n$ can be replaced by its natural limit $\eta_0$ uniformly in some other indexing set $\Theta$. In particular we reconsider some examples treated by Ghoudi and Remillard [*Asymptotic Methods in Probability and Statistics* (1998) 171–197, *Fields Inst. Commun.* 44 (2004) 381–406]. We recast their examples in terms of empirical process theory, and provide an alternative general view which should be of wide applicability.

## 1. Introduction

Let $X_1, \ldots, X_n$ be i.i.d. random elements in a measurable space $(\mathcal{X}, \mathcal{A})$ with law $P$, and for a measurable function $f : \mathcal{X} \to \mathbb{R}$ let the expectation, empirical measure and empirical process at $f$ be denoted by

$$Pf = \int f \, dP, \qquad \mathbb{P}_n f = \frac{1}{n} \sum_{i=1}^{n} f(X_i), \qquad \mathbb{G}_n f = \sqrt{n}(\mathbb{P}_n - P)f.$$

Given a collection $\{f_{\theta,\eta} : \theta \in \Theta, \eta \in H\}$ of measurable functions $f_{\theta,\eta} : \mathcal{X} \to \mathbb{R}$ indexed by sets $\Theta$ and $H$ and "estimators" $\eta_n$, we wish to prove that, as $n \to \infty$,

$$(1) \qquad \sup_{\theta \in \Theta} \left| \mathbb{G}_n(f_{\theta,\eta_n} - f_{\theta,\eta_0}) \right| \to_p 0.$$

Here an "estimator" $\eta_n$ is a random element with values in $H$ defined on the same probability space as $X_1, \ldots, X_n$, and $\eta_0 \in H$ is a fixed element, which is typically a limit in probability of the sequence $\eta_n$.

The result (1) is interesting for several applications. A direct application is to the estimation of the functional $\theta \mapsto Pf_{\theta,\eta}$. If the parameter $\eta$ is unknown, we may replace it by an estimator $\eta_n$ and use the empirical estimator $\mathbb{P}_n f_{\theta,\eta_n}$. The result (1) helps to derive the limit behaviour of this estimator, as we can decompose

$$(2) \qquad \sqrt{n}(\mathbb{P}_n f_{\theta,\eta_n} - Pf_{\theta,\eta_0}) = \mathbb{G}_n(f_{\theta,\eta_n} - f_{\theta,\eta_0}) + \mathbb{G}_n f_{\theta,\eta_0} + \sqrt{n}P(f_{\theta,\eta_n} - f_{\theta,\eta_0}).$$

If (1) holds, then the first term on the right converges to zero in probability. Under appropriate conditions on the functions $f_{\theta,\eta_0}$, the second term on the right

*Supported in part by NSF Grant DMS-05-03822, NI-AID Grant 2R01 AI291968-04, and by grant B62-596 of the Netherlands Organisation of Scientific Research NWO

[1] Section Stochastics, Department of Mathematics, Faculty of Sciences, Vrije Universiteit, De Boelelaan 1081a, 1081 HV Amsterdam, e-mail: `aad@cs.vu.nl`

[2] University of Washington, Department of Statistics, Box 354322, Seattle, Washington 98195-4322, USA, e-mail: `jaw@stat.washington.edu`

*AMS 2000 subject classifications:* 62G07, 62G08, 62G20, 62F05, 62F15.

*Keywords and phrases:* delta-method, Donsker class, entropy integral, pseudo observation.

will converge to a Gaussian process by the (functional) central limit theorem. The behavior of the third term depends on the estimators $\eta_n$, and would typically follow from an application of the (functional) delta-method, applied to the map $\eta \mapsto (Pf_{\theta,\eta} : \theta \in \Theta)$.

In an interesting particular case of this situation, the functions $f_{\theta,\eta}$ take the form

$$f_{\theta,\eta}(x) = \theta\big(\eta(x)\big),$$

for maps $\theta : \mathbb{R}^d \to \mathbb{R}$ and each $\eta \in H$ being a map $\eta : \mathcal{X} \to \mathbb{R}^d$. The realizations of the estimators $\eta_n$ are then functions $x \mapsto \eta_n(x) = \eta_n(x; X_1, \ldots, X_n)$ on the sample space $\mathcal{X}$ and can be evaluated at the observations to obtain the random vectors $\eta_n(X_1), \ldots, \eta_n(X_n)$ in $\mathbb{R}^d$. The process $\{\mathbb{P}_n f_{\theta,\eta_n} : \theta \in \Theta\}$ is the empirical measure of these vectors indexed by the functions $\theta$. For instance, if $\Theta$ consists of the indicator functions $1_{(-\infty,\theta]}$ for $\theta \in \mathbb{R}^d$, then this measure is the empirical distribution function

$$\theta \mapsto \mathbb{P}_n f_{\theta,\eta_n} = \frac{1}{n}\sum_{i=1}^n 1\{\eta_n(X_i) \leq \theta\}$$

of the random vectors $\eta_n(X_1), \ldots, \eta_n(X_n)$. The properties of such empirical processes were studied in some generality and for examples of particular interest in Ghoudi and Remillard [6, 7]. Ghoudi and Remillard [6] apparently coined the name "pseudo-observations" for the vectors $\eta_n(X_1), \ldots, \eta_n(X_n)$. The examples include, for instance, regression residuals, Kendall's dependence process, and copula processes; see the end of Section 2 for explicit formulation of these three particular examples. One purpose of the present paper is to extend the results in these papers also to other index classes $\Theta$ besides the class of indicator functions. Another purpose is to recast their results in terms of empirical process theory, which leads to simplification and alternative conditions.

A different, indirect application of (1) is to the derivation of the asymptotic distribution of $Z$-estimators. A $Z$-estimator for $\theta$ might be defined as the solution $\hat{\theta}_n$ of the equation $\mathbb{P}_n f_{\theta,\eta_n} = 0$, where again an unknown "nuisance" parameter $\eta$ is replaced by an estimator $\eta_n$. In this case (1) shows that

$$\mathbb{P}_n f_{\hat{\theta}_n,\eta_n} - \mathbb{P}_n f_{\hat{\theta}_n,\eta_0} = P(f_{\hat{\theta}_n,\eta_n} - f_{\hat{\theta}_n,\eta_0}) + o_P(1/\sqrt{n}),$$

so that the limit behavior of $\hat{\theta}_n$ can be derived by comparison with the estimating equation defined by $\mathbb{P}_n f_{\theta,\eta_0}$ (with $\eta_0$ substituted for $\eta_n$). The "drift" sequence $P(f_{\hat{\theta}_n,\eta_n} - f_{\hat{\theta}_n,\eta_0})$, which will typically be equivalent to $P(f_{\theta_0,\eta_n} - f_{\theta_0,\eta_0})$ up to order $o_P(1/\sqrt{n})$, may give rise to an additional component in the limit distribution.

The paper is organized as follows. In Section 2 we derive general conditions for the validity of (1) and formulate several particular examples to be considered in more detail in the sequel. In Section 3 we specialize the general results to composition maps. In Section 4 we combine these results with results on Hadamard differentiability to obtain the asymptotic distribution of empirical processes indexed by pseudo observations. Finally in Section 5 we formulate our results for several of the particular examples mentioned above and at the end of Section 2.

## 2. General result

In many situations we wish to establish (1) without knowing much about the nature of the estimators $\eta_n$, beyond possibly that they are consistent for some value $\eta_0$.

For instance, this is true if (1) is used as a step in the derivation of $M-$ or $Z-$ estimators. (Cf. Van der Vaart and Wellner [12] and Van der Vaart [11].) Then an appropriate method of establishing (1) is through a Donsker or entropy condition, as in the following theorems. Proofs of the Theorems 2.1 and 2.2 can be found in the mentioned references.

Both theorems assume that $\eta_n$ is "consistent for $\eta_0$" in the sense that

$$(3) \qquad \sup_{\theta \in \Theta} P(f_{\theta,\eta_n} - f_{\theta,\eta_0})^2 \to_p 0.$$

**Theorem 2.1.** *Suppose that $H_0$ is a fixed subset of $H$ such that $\Pr(\eta_n \in H_0) \to 1$ and suppose that the class of functions $\{f_{\theta,\eta} : \theta \in \Theta, \eta \in H_0\}$ is $P$-Donsker. If (3) holds, then (1) is valid.*

For the second theorem, let $N(\epsilon, \mathcal{F}, L_2(P))$ and $N_{[\,]}(\epsilon, \mathcal{F}, L_2(P))$ be the $\epsilon$-covering and $\epsilon$-bracketing numbers of a class $\mathcal{F}$ of measurable functions (cf. Pollard [8] and van der Vaart and Wellner [12]) and define entropy integrals by

$$(4) \qquad J(\delta, \mathcal{F}, L_2) = \int_0^\delta \sup_Q \sqrt{\log N(\epsilon \|F\|_{Q,2}, \mathcal{F}, L_2(Q))} \, d\epsilon,$$

$$(5) \qquad J_{[\,]}(\delta, \mathcal{F}, L_2(P)) = \int_0^\delta \sqrt{\log N_{[\,]}(\epsilon \|F\|_{P,2}, \mathcal{F}, L_2(P))} \, d\epsilon.$$

Here $F$ is an arbitrary, measurable envelope function for the class $\mathcal{F}$: a measurable function $F : \mathcal{X} \to \mathbb{R}$ such that $|f(x)| \le F(x)$ for every $f \in \mathcal{F}$ and $x \in \mathcal{X}$. We say that a sequence $F_n$ of envelope functions satisfies the Lindeberg condition if $PF_n^2 = O(1)$ and $PF_n^2 1_{F_n \ge \epsilon \sqrt{n}} \to 0$ for every $\epsilon > 0$.

**Theorem 2.2.** *Suppose that $H_n$ are subsets of $H$ such that $\Pr(\eta_n \in H_n) \to 1$ and such that the classes of functions $\mathcal{F}_n = \{f_{\theta,\eta} : \theta \in \Theta, \eta \in H_n\}$ satisfy either $J_{[\,]}(\delta_n, \mathcal{F}_n, L_2(P)) \to 0$, or $J(\delta_n, \mathcal{F}_n, L_2) \to 0$ for every sequence $\delta_n \to 0$, relative to envelope functions that satisfy the Lindeberg condition. In the second case also assume that the classes $\mathcal{F}_n$ are suitably measurable (e.g. countable). If (3) holds, then (1) is valid.*

Because there are many techniques to verify that a given class of functions is Donsker, or to compute bounds on its entropy integrals, the preceding lemmas give quick results, if they apply. Furthermore, they appear to be close to best possible unless more information about the estimators $\eta_n$ can be brought in, or explicit computations are possible for the functions $f_{\theta,\eta}$.

In some applications the estimators $\eta_n$ are known to converge at a certain rate and/or known to possess certain regularity properties (e.g. uniform bounded derivatives). Such knowledge cannot be exploited in Theorem 2.1, but could be used for the choice of the sets $H_n$ in Theorem 2.2. We now discuss an alternative approach which can be used if the estimators $\eta_n$ are also known to converge in distribution, if properly rescaled.

Let $H$ be a Banach space, and suppose that the sequence $\sqrt{n}(\eta_n - \eta_0)$ converges in distribution to a tight, Borel-measurable random element in $H$. The "convergence in distribution" may be understood in the sense of Hoffmann-Jørgensen, so that $\eta_n$ need not be Borel-measurable itself.

The tight limit of the sequence $\sqrt{n}(\eta_n - \eta_0)$ takes its values in a $\sigma$-compact subset $H_0 \subset H$. For $\theta \in \Theta$, $h_0 \in H_0$, and $\delta > 0$ define a sequence of classes of functions by

$$(6) \qquad \mathcal{F}_n(\theta, h_0, \delta) = \{f_{\theta,\eta_0 + n^{-1/2}h} - f_{\theta,\eta_0 + n^{-1/2}h_0} : h \in H, \ \|h - h_0\| < \delta\}.$$

Let $F_n(\theta, h_0, \delta)$ be arbitrary measurable envelope functions for these classes.

**Theorem 2.3.** *Suppose that the sequence $\sqrt{n}(\eta_n - \eta_0)$ converges in distribution to a tight, random element with values in a given $\sigma$-compact subset $H_0$ of $H$. Suppose that*

  (i) $\sup_\theta |\mathbb{G}_n(f_{\theta, \eta_0 + n^{-1/2} h_0} - f_{\theta, \eta_0})| \to_p 0$ *for every $h_0 \in H_0$.*
  (ii) $\sup_\theta |\mathbb{G}_n F_n(\theta, h_0, \delta)| \to_p 0$ *for every $\delta > 0$ and every $h_0 \in H_0$;*
  (iii) $\sup_\theta \sup_{h_0 \in K} \sqrt{n} \, P F_n(\theta, h_0, \delta_n) \to 0$ *for every $\delta_n \to 0$ and every compact $K \subset H_0$;*

*Then* (1) *is valid.*

*Proof.* Suppose that $\sqrt{n}(\eta_n - \eta_0) \Rightarrow Z$ and let $\epsilon > 0$ be fixed. There exists a compact set $K \subset H_0$ with $P(Z \in K) > 1 - \epsilon$ and hence for every $\delta > 0$, with $K^\delta$ the set of all points at distance less than $\delta$ to $K$,

$$\liminf_{n \to \infty} \Pr\big(\sqrt{n}(\eta_n - \eta_0) \in K^{\delta/2}\big) > 1 - \epsilon.$$

In view of the compactness of $K$ there exist finitely many elements $h_1, \ldots, h_p \in K \subset H_0$ (with $p = p(\delta)$ depending on $\delta$) such that the balls of radius $\delta/2$ around these points cover $K$. Then $K^{\delta/2}$ is contained in the union of the balls of radius $\delta$, by the triangle inequality. Thus, with $B(h, \delta)$ denoting the ball of radius $\delta$ around $h$ in the space $H$,

$$\big\{\sqrt{n}(\eta_n - \eta_0) \in K^{\delta/2}\big\} \subset \bigcup_{i=1}^{p(\delta)} \big\{\eta_n \in B(\eta_0 + n^{-1/2} h_i, \delta)\big\}.$$

It follows that with probability at least $1 - \epsilon$, as $n \to \infty$,

$$\sup_\theta |\mathbb{G}_n(f_{\theta, \eta_n} - f_{\theta, \eta_0})|$$
$$\leq \sup_\theta \max_i \sup_{\|h - h_i\| < \delta} |\mathbb{G}_n(f_{\theta, \eta_0 + n^{-1/2} h} - f_{\theta, \eta_0})|$$
$$\leq \sup_\theta \max_i \sup_{\|h - h_i\| < \delta} \Big[ |\mathbb{G}_n(f_{\theta, \eta_0 + n^{-1/2} h} - f_{\theta, \eta_0 + n^{-1/2} h_i})| $$
$$+ |\mathbb{G}_n(f_{\theta, \eta_0 + n^{-1/2} h_i} - f_{\theta, \eta_0})| \Big]$$
$$\leq \sup_\theta \max_i |\mathbb{G}_n F_n(\theta, h_i, \delta)| + 2 \sup_\theta \sup_{h_0 \in K} \sqrt{n} \, P F_n(\theta, h_0, \delta)$$
$$+ \sup_\theta \max_i |\mathbb{G}_n(f_{\theta, \eta_0 + n^{-1/2} h_i} - f_{\theta, \eta_0})|,$$

where in the last step we use the inequality $|\mathbb{G}_n f| \leq |\mathbb{G}_n F| + 2\sqrt{n} PF$, valid for any functions $f$ and $F$ with $|f| \leq F$. The maxima in the display are over the finite set $i = 1, \ldots, p(\delta)$, and the elements $h_1, \ldots, h_{p(\delta)} \in K$ depend on $\delta$. By assumptions (i) and (ii) the first and third terms converge to zero as $n \to \infty$, for every fixed $\delta$. It follows that there exists $\delta_n \downarrow 0$ such that these terms with $\delta_n$ substituted for $\delta$ converge to 0. For this $\delta_n$, all three terms converge to zero in probability as $n \to \infty$.  $\square$

The rate of convergence $\sqrt{n}$ in the preceding theorem may be replaced by another rate, with appropriate changes in the conditions, but the rate $\sqrt{n}$ appears natural in the following context. For more general metrizable topological vector spaces there are similar, but less attractive, results possible.

The two conditions (i), (ii) of Theorem 2.3 concern the empirical process indexed by the classes of functions

$$\{f_{\theta, \eta_0 + n^{-1/2} h_0} - f_{\theta, \eta_0} : \theta \in \Theta\}, \tag{7}$$

$$\{F_n(\theta, h_0, \delta) : \theta \in \Theta\}. \tag{8}$$

These classes are indexed by $\Theta$ only, and hence Theorem 2.3, if applicable, avoids conditions for (1) that involve measures of the complexity of the class $\{f_{\theta, \eta} : \theta \in \Theta, \eta \in H\}$ due to the parameter $\eta \in H$.

Condition (iii) of Theorem 2.3 involves the mean of the envelopes of the classes $\mathcal{F}_n(\theta, h_0, \delta)$. For the minimal envelopes this condition takes the form

$$\sup_{\theta} \sup_{h_0 \in K} \sqrt{n} \, P \sup_{\|h - h_0\| < \delta_n} |f_{\theta, \eta_0 + n^{-1/2} h} - f_{\theta, \eta_0 + n^{-1/2} h_0}| \to 0 \tag{9}$$

for all $\delta_n \downarrow 0$. This is an "integrated uniform local Lipschitz assumption" on the dependence $\eta \mapsto f_{\theta, \eta}$. In some applications it may be useful not to use the minimal envelope functions. The lemma is valid for any envelope functions, as long as the same envelopes are used in both (ii) and (iii).

The set $K$ in (iii) or (9) is a compact set in the support of the limit distribution of the sequence $\sqrt{n}(\eta_n - \eta_0)$. In some cases condition (iii) may be valid for any compact $K \subset H$, whereas in other cases more precise information about the limit process must be exploited. For instance, if the sequence $\sqrt{n}(\eta_n - \eta_0)$ converges in distribution to a tight zero-mean Gaussian process $G$ in the space $H = \ell^\infty(T)$ of bounded functions on some set $T$, then $K$ may be taken to be a set of functions $z : T \to \mathbb{R}$ that is uniformly bounded and uniformly equicontinuous relative to the semimetric with square $d^2(s, t) = \mathrm{E}(G_s - G_t)^2$ (and $T$ will be totally bounded for $d$). Cf. e.g. van der Vaart and Wellner [12], page 39.

Condition (iii) is an analytical condition, whereas conditions (i) and (ii) are empirical process conditions. In many cases the latter pair of conditions can be verified by standard empirical process type arguments. For reference we quote two lemmas that allow handling the empirical process indexed by a sequence of classes, as in (8) or (7). (For proofs see e.g. van der Vaart [10, 11].) Both lemmas apply to classes $\mathcal{F}_n$ of measurable functions $f : \mathcal{X} \mapsto \mathbb{R}$ such that

$$\sup_{f \in \mathcal{F}_n} P f^2 \to 0. \tag{10}$$

**Lemma 2.1.** *Suppose that the class of functions $\bigcup_n \mathcal{F}_n$ is $P$-Donsker. If (10) holds, then $\sup_{f \in \mathcal{F}_n} |\mathbb{G}_n(f)| \to_p 0$.*

**Lemma 2.2.** *Suppose that either $J_{[\cdot]}(\delta_n, \mathcal{F}_n, L_2(P)) \to 0$ or $J(\delta_n, \mathcal{F}_n, L_2) \to 0$ for all $\delta_n \downarrow 0$ relative to envelope functions $F_n$ that satisfy the Lindeberg condition. In the second case also assume that each class $\mathcal{F}_n$ is suitably measurable. If (10) holds, then $\sup_{f \in \mathcal{F}_n} |\mathbb{G}_n(f)| \to_p 0$.*

**Example 1** (Regression residual processes). Suppose that $(X_1, Y_1), \ldots, (X_n, Y_n)$ are a random sample distributed according to the regression model $Y = g_\eta(X) + e$. For given estimators $\eta_n$ we can form the residuals $\hat{e}_i = Y_i - g_{\eta_n}(X_i)$ and may be interested in the empirical process corresponding to $\hat{e}_1, \ldots, \hat{e}_n$, i.e. for a collection $\Theta$ of functions $\theta : \mathbb{R} \to \mathbb{R}$ we consider the process $\{n^{-1} \sum_{i=1}^n \theta(\hat{e}_i) : \theta \in \Theta\}$. This fits the general set-up with the functions $f_{\theta, \eta}$ defined as $f_{\theta, \eta}(x, y) = \theta(y - g_\eta(x))$.

In many cases it will be possible to apply Theorem 2.1. For instance, if $x \in \mathbb{R}^d$, $g_\eta$ is a polynomial in $x$, and $\Theta$ is the class of indicator functions $1_{(-\infty, \theta]}$ for $\theta \in \mathbb{R}$, then the functions $f_{\theta, \eta}$ are the indicator functions of the sets $\{(x, y) : y - g_\eta(x) - \theta \leq 0\}$. Because the set of functions $(x, y) \mapsto y - g_\eta(x) - \theta$ is contained in a finite-dimensional vector space, it is a VC-class, and hence so are their negativity sets (e.g. van der Vaart and Wellner [12], Lemma 2.6.18). Thus the class of functions $f_{\theta, \eta}$ is Donsker, and Theorem 2.1 can be applied directly.

**Example 2** (Kendall's process). Let $\eta_n$ be the empirical distribution function of a random sample $X_1, \ldots, X_n$ from a distribution $\eta_0$ on $\mathbb{R}^d$. Barbe, Genest, Ghoudi and Remillard [2] and Ghoudi and Remillard [6] study the behavior of the empirical distribution function $K_n$ of the pseudo-observations $\eta_n(X_i)$,

$$K_n(\theta) = \frac{1}{n} \sum_{i=1}^{n} 1\{\eta_n(X_i) \leq \theta\}, \qquad \theta \in [0, 1],$$

and the resulting Kendall's process

$$(11) \qquad \qquad \sqrt{n}(K_n(\theta) - K(\theta)), \qquad \theta \in [0, 1]$$

where $K(\theta) = P(\eta_0(X) \leq \theta)$. This fits the general set-up with $f_{\theta, \eta}$ the composition function $f_{\theta, \eta} = \theta \circ \eta$, and $\theta$ the indicator function $1_{(-\infty, \theta]}$ (where we abuse notation by using the symbol $\theta$ in two different ways).

An attempt to apply Theorem 2.1 to this problem would lead to the consideration of the class of all indicator functions of sets of the form $\{x \in \mathbb{R}^d : \eta(x) \leq \theta\}$ for $\eta$ ranging over the cumulative distribution functions on $\mathbb{R}^d$ and $\theta \in [0, 1]$. This class is similar to the collection of all "lower layers" in $\mathbb{R}^d$, and, unfortunately, fails to be Donsker for most distributions (cf. Dudley [3], page 264, 373 or Dudley [4]). In this case it appears to be necessary to exploit the limit behaviour of the sequence $\sqrt{n}(\eta_n - \eta_0)$. Ghoudi and Remillard [6] have shown that (1) is valid in this case, under some strong smoothness assumptions on the underlying measure $\eta_0$. In Sections 4 and 5 we rederive some of their results by empirical process methods using Theorem 2.3.

We also consider the empirical process of the variables $\eta_n(X_1), \ldots, \eta_n(X_n)$ indexed by classes of functions other than the indicators $1_{(-\infty, \theta]}$. If the indexing functions are smooth, then this empirical process will converge even without smoothness conditions on $\eta_0$. A proof can be based on Theorem 2.3.

**Example 3** (Copula processes). Suppose that $X_1, \ldots, X_n$ are a sample from a distribution $\eta_0$ on $\mathbb{R}^d$. Write $X_i = (X_{i,1}, \ldots, X_{i,d})$ and let $\eta_{0,1}, \ldots, \eta_{0,d}$ be the marginal distributions. The *copula* function $C$ associated with $\eta_0$ is the distribution function of the vector $(\eta_{0,1}(X_{1,1}), \ldots, \eta_{0,d}(X_{1,d}))$, i.e. with $\eta_{0,j}^{-1}(u) = \inf\{x : \eta_{0,j}(x) \geq u\}$ for $u \in [0, 1]$,

$$C(u_1, \ldots, u_d) = \eta_0(\eta_{0,1}^{-1}(u_1), \ldots, \eta_{0,d}^{-1}(u_d))$$

for $(u_1, \ldots, u_d) \in [0, 1]^d$. For $j = 1, \ldots, d$ let $\eta_{n,j}$ be the empirical distribution function of $X_{1,j}, \ldots, X_{n,j}$ (on $\mathbb{R}$), and let $\eta_n$ be the empirical distribution function of $X_1, \ldots, X_n$ (on $\mathbb{R}^d$). Then a natural estimator $C_n$ of $C$ is given by

$$C_n(u) = \frac{1}{n} \sum_{i=1}^{n} 1\{\hat{e}_{n,i} \leq u\}, \qquad u \in [0, 1]^d,$$

for the "pseudo-observations" $\hat{e}_{n,i} = \big(\eta_{n,1}(X_{i,1}), \ldots, \eta_{n,d}(X_{i,d})\big)$. The resulting "copula processes",

$$(12) \qquad \sqrt{n}(C_n(u) - C(u)), \qquad u \in [0,1]^d,$$

have been considered by Stute [9], Gänssler and Stute [5], and Ghoudi and Remillard [7]. This example can be treated using Theorem 2.3, but also with the more straightforward Theorem 2.1, or even by employing the theory of Hadamard differentiability, as in Chapter 3.9 of Van der Vaart and Wellner [12].

## 3. Composition

In this section we consider the case where the functions $f_{\theta,\eta}$ take the form

$$(13) \qquad f_{\theta,\eta}(x) = \theta(\eta(x)),$$

for $\theta$ ranging over a class $\Theta$ of functions $\theta : \mathbb{R}^d \to \mathbb{R}$ and $\eta$ ranging over a class $H$ of measurable functions $\eta : \mathcal{X} \to \mathbb{R}^d$. We first give general conditions for the validity of condition (iii) of Theorem 2.3, and next consider also the conditions (i) and (ii) for the special cases of functions $\theta$ that are Lipschitz and monotone, respectively. We develop these results for the case that the sequence $\sqrt{n}(\eta_n - \eta_0)$ converges in distribution in the space $H = \ell^\infty(\mathcal{X}, \mathbb{R}^d)$ of uniformly bounded functions $z : \mathcal{X} \to \mathbb{R}^d$, equipped with the uniform norm $\|z\| = \sup_{x \in \mathcal{X}} \|z(x)\|$. (Variations of these results are possible. For instance, $\mathbb{R}^d$ could be replaced by a more general Banach space, and $H$ could be equipped with a weighted uniform norm.)

### 3.1. *Condition* (i)

For $f_{\theta,\eta}$ taking the form (13), condition (i) of Theorem 2.3 takes the form

$$(14) \qquad \sup_{f \in \mathcal{F}_n} |\mathbb{G}_{n,Q}f| \to_Q 0$$

for $Q = P \circ (\eta_0, h_0)^{-1}$, $\mathbb{G}_{n,Q}$ the empirical process of a random sample from the measure $Q$, and $\mathcal{F}_n$ the class of functions

$$(15) \qquad \mathcal{F}_n = \big\{(y,z) \mapsto \theta(y + n^{-1/2}z) - \theta(y) : \theta \in \Theta\big\}.$$

Condition (i) requires that (14) is valid for every fixed choice of $h_0 \in H_0$, i.e. for every measure $Q$ determined as the law of $(\eta_0(X), h_0(X))$ for some $h_0 \in H_0$ and $X$ distributed according to $P$.

This situation is of the form considered in Lemmas 2.1 and 2.2, and both lemmas may be applicable in a given setting. It is not especially helpful to restate these lemmas for the present special situation. Instead, we give one easy to check set of sufficient conditions. This covers VC-classes $\Theta$, and much more.

If $\Theta_{env} : \mathbb{R}^d \to \mathbb{R}$ is an envelope function for $\Theta$, then

$$(16) \qquad F_n(y,z) = \Theta_{env}(y + n^{-1/2}z) + \Theta_{env}(y)$$

is an envelope function for $\mathcal{F}_n$. (A crude one, because we do not exploit that the functions in $\mathcal{F}_n$ are differences.)

**Lemma 3.1.** *Suppose that $J(1, \Theta, L_2) < \infty$, that $\Theta$ is suitably measurable, that $P(\Theta_{env} \circ \eta_0)^2 < \infty$, and that the functions $\Theta_{env} \circ (\eta_0 + n^{-1/2} h_0)$ satisfy the Lindeberg condition in $L_2(P)$, for every $h_0 \in H_0$. If*

$$\sup_{\theta \in \Theta} P\big(\theta \circ (\eta_0 + n^{-1/2} h_0) - \theta \circ \eta_0\big)^2 \to 0,$$

*for every $h_0 \in H_0$, then condition (i) of Theorem 2.3 is satisfied for the functions $f_{\theta, \eta}$ given by (13).*

*Proof.* It suffices to prove (14) for the classes $\mathcal{F}_n$ given in (15). The class $\mathcal{F}_n$ is contained in the difference of the classes

$$\mathcal{F}_n' = \{(y, z) \mapsto \theta(y + n^{-1/2} z) : \theta \in \Theta\},$$
$$\mathcal{F}'' = \{(y, z) \mapsto \theta(y) : \theta \in \Theta\}.$$

These classes possess envelope functions $F_n'$ and $F''$ defined by

$$F_n'(y, z) = \Theta_{env}(y + n^{-1/2} z),$$
$$F''(y, z) = \Theta_{env}(y).$$

The uniform entropy of $\mathcal{F}''$ relative to $F''$ is finite by assumption. The uniform entropy of $\mathcal{F}_n'$ relative to $F_n'$ is exactly the same, as the law of $Y + n^{-1/2} Z$ runs through all possible laws on $\mathbb{R}^d$ if the law of $(Y, Z)$ runs through all possible laws on $\mathbb{R}^d \times \mathbb{R}^d$. The uniform entropy of $\mathcal{F}_n$ relative to $F_n$ is bounded by the sum of the uniform entropies of $\mathcal{F}_n'$ and $\mathcal{F}''$. (Cf. e.g. Theorem 2.10.20 of van der Vaart and Wellner [12].) Now apply Lemma 2.2. $\qquad \square$

### 3.2. Lipschitz functions $\theta$

Assume that every function $\theta : \mathbb{R}^d \to \mathbb{R}$ in the class $\Theta$ is uniformly Lipschitz in that

$$(17) \qquad |\theta(r_1) - \theta(r_2)| \leq \|r_1 - r_2\|.$$

Then, for every $x \in \mathcal{X}$,

$$\left| \theta\big(\eta_0(x) + n^{-1/2} h(x)\big) - \theta\big(\eta_0(x) + n^{-1/2} h_0(x)\big) \right| \leq \frac{\|h(x) - h_0(x)\|}{\sqrt{n}}.$$

The norm in the right side is bounded by the supremum norm $\|h - h_0\|$ on $\ell^\infty(\mathcal{X}, \mathbb{R}^d)$. It follows that the classes $\mathcal{F}_n(\theta, h_0, \delta)$ as in (6) possess envelope functions

$$(18) \qquad F_n(\theta, h_0, \delta) = \delta / \sqrt{n}.$$

**Theorem 3.1.** *If $\Theta$ is a suitably measurable collection of uniformly bounded, uniformly Lipschitz functions $\theta : \mathbb{R}^d \to \mathbb{R}$ such that $J(1, \Theta, L_2) < \infty$ (relative to a constant envelope function), $\eta_0 \in \ell^\infty(\mathcal{X}, \mathbb{R}^d)$, and the sequence $\sqrt{n}(\eta_n - \eta_0)$ converges weakly in $\ell^\infty(\mathcal{X}, \mathbb{R}^d)$ to a tight random element, then*

$$\sup_{\theta \in \Theta} \big| \mathbb{G}_n\big(\theta(\eta_n) - \theta(\eta_0)\big) \big| \to_p 0.$$

*Proof.* With the envelope functions $F_n(\theta, h_0, \delta)$ as defined in (18), condition (ii) of Theorem 2.3 is trivially satisfied because the envelopes are actually constants, and the validity of condition (iii) is immediate.

By assumption we can choose the envelope function of $\Theta$ equal to a constant and $J(1, \Theta, L_2) < \infty$. This suffices for the verification of most of the conditions of Lemma 3.1. Finally, it suffices to note that

$$P\big(\theta \circ (\eta_0 + n^{-1/2}h_0) - \theta \circ \eta_0\big)^2 \le \|h_0\|^2/n.$$

By Lemma 3.1 we conclude that condition (i) of Theorem 2.3 is also satisfied, whence the theorem follows from Theorem 2.3. □

For the verification of condition (i) of Theorem 2.3 it suffices to consider the functions $\theta$ on the range of the functions $\eta_0 + h_0/\sqrt{n}$ for a fixed $h_0$ in the support of the limit distribution of the sequence $\sqrt{n}(\eta_n - \eta_0)$. Thus we may restrict the functions $\theta$ to a subset of $\mathbb{R}^d$ that contains the ranges of these functions and interpret the condition $J(1, \Theta, L_2) < \infty$ in Lemma 3.1 accordingly. In particular, in Theorem 2.3 we may replace this condition by the condition that $J(1, \Theta_K, L_2) < \infty$ for every norm-bounded subset $K \subset \mathbb{R}^d$, where $\Theta_K$ is the collection of restrictions $\theta : K \to \mathbb{R}$ of the functions $\theta \in \Theta$.

Any collection of uniformly bounded, Lipschitz functions $\theta : K \to \mathbb{R}$ on a compact interval $K$ satisfies $J(1, \Theta, L_2) < \infty$. (Cf. e.g. van der Vaart and Wellner [12], page 157.) Thus in the case that $d = 1$ the assertion of the preceding theorem is true for any collection of uniformly bounded Lipschitz functions.

For $d > 1$ further restrictions on the class $\Theta$ may be necessary. For instance, any subset of the unit ball in the Hölder space $C^\alpha(K)$ for a compact interval $K \subset \mathbb{R}^d$ possesses a finite uniform entropy integral provided $\alpha > d/2$. (Cf. e.g. van der Vaart and Wellner [12], page 157.) The assertion of the preceding theorem is also true for such a class.

There are many other examples of classes of Lipschitz functions with finite uniform entropy integrals. For instance, VC-classes of Lipschitz functions.

### 3.3. Monotone functions $\theta$

Assume that every function $\theta : \mathbb{R}^d \to \mathbb{R}$ in $\Theta$ is the survival function $\theta(x) = \int 1_{[x,\infty)} d\theta$ of a subprobability measure on $\mathbb{R}^d$. Then each $\theta$ is nonincreasing in each of its arguments. If $H = \ell^\infty(\mathcal{X}, \mathbb{R}^d)$ equipped with the uniform norm relative to the max-norm on $\mathbb{R}^d$, then

$$\Big|\theta(\eta_0 + n^{-1/2}h) - \theta(\eta_0 + n^{-1/2}h_0)\Big|$$
$$\le \theta(\eta_0 + n^{-1/2}h_0 - n^{-1/2}\|h - h_0\|) - \theta(\eta_0 + n^{-1/2}h_0 + n^{-1/2}\|h - h_0\|).$$

It follows that the classes $\mathcal{F}_n(\theta, h_0, \delta)$ possess envelope functions, with $\delta$ the vector $(\delta, \dots, \delta)$,

$$F_n(\theta, h_0, \delta) = \theta(\eta_0 + n^{-1/2}h_0 - n^{-1/2}\delta) - \theta(\eta_0 + n^{-1/2}h_0 + n^{-1/2}\delta).$$

In order to verify condition (iii) of Theorem 2.3, we assume that for given (possibly infinite) $a < b$ in $\mathbb{R}^d$ and every $\delta_n \downarrow 0$ and compact set $K \subset H_0 \cup \{0\}$,

$$(19) \qquad \sup_{t \in \mathbb{R}^d, a \le t \le b} \sup_{h_0 \in K} \sqrt{n} P\big(1_{\eta_0 + n^{-1/2}h_0 \le t + n^{-1/2}\delta_n} - 1_{\eta_0 + n^{-1/2}h_0 \le t}\big) \to 0.$$

**Theorem 3.2.** *Let $\Theta$ be a collection of survival functions $\theta : \mathbb{R}^d \to [0,1]$ of subprobability measures supported on an interval $(a,b) \subset \mathbb{R}^d$. If the sequence $\sqrt{n}(\eta_n - \eta_0)$ converges in distribution in $\ell^\infty(\mathcal{X}, \mathbb{R}^d)$ to a tight Borel measure concentrating on the $\sigma$-compact set $H_0$, and (19) holds for every $\delta_n \downarrow 0$ and every compact $K \subset H_0 \cup \{0\}$, then*

$$\sup_\theta \big| \mathbb{G}_n \big( \theta(\eta_n) - \theta(\eta_0) \big) \big| \to_p 0.$$

*Proof.* The survival functions of subprobability measures are in the convex hull of the set of indicator functions $1_{[t,\infty)}$, which is a VC-class. Therefore the entropy integral $J(1, \Theta, L_2)$ of $\Theta$ relative to a constant envelope is finite. (Cf. e.g. van der Vaart and Wellner [12], page 145.)

Defining $F_n(\theta, \eta_0, \delta)$ as in the display preceding the theorem, we can write

$$PF_n(\theta, \eta_0, \delta)$$
$$= \int P\big(1_{(-\infty, s]}(\eta_0 + n^{-1/2}h_0 - n^{-1/2}\delta) - 1_{(-\infty, s]}(\eta_0 + n^{-1/2}h_0)\big)\, d\theta(s)$$
$$\leq \|\theta\| \sup_{a \leq s \leq b} P\big(1_{(-\infty, s]}(\eta_0 + n^{-1/2}h_0 - n^{-1/2}\delta) - 1_{(-\infty, s]}(\eta_0 + n^{-1/2}h_0)\big).$$

By assumption (19) the right side converges to zero faster than $1/\sqrt{n}$, for every $\delta = \delta_n \downarrow 0$, uniformly in $h_0 \in K$, and uniformly in $\theta$ because the total variation norms $\|\theta\|$ are uniformly bounded. This verifies condition (iii) of Theorem 2.3.

Because $\theta$ is monotone with range contained in $[0,1]$,

$$P\big(\theta(\eta_0 + n^{-1/2}\delta) - \theta(\eta_0)\big)^2 \leq \sup_{a \leq s \leq b} P\big(1_{(-\infty, s]}(\eta_0 - n^{-1/2}\delta) - 1_{(-\infty, s]}(\eta_0)\big).$$

By assumption (19) with $h_0 = 0$ this converges to zero faster than $1/\sqrt{n}$ for every sequence $\delta_n \downarrow 0$. This can be seen to imply that the expression in the display (which does not have the leading $\sqrt{n}$) converges to zero also for fixed $\delta$. By monotonicity of $\theta$ we can bound $\big|\theta(\eta_0 + h_0/\sqrt{n}) - \theta(\eta_0)\big|$ by $\big|\theta(\eta_0 - \delta/\sqrt{n}) - \theta(\eta_0)\big|$ for $\delta = \|h_0\|$. By Lemma 3.1 we now conclude that condition (i) of Theorem 2.3 is satisfied.

In the present case the envelope functions $F_n(\theta, h_0, \delta)$ are equal to the functions $f_{\theta, \eta_0 + n^{-1/2}h} - f_{\theta, \eta_0 - n^{-1/2}h}$ for $h = h_0 - \delta$. Therefore, the validity of condition (ii) of Theorem 2.3 follows by the same arguments as used for the validity of condition (i). □

Condition (19) is a uniform Lipschitz condition on the distribution functions of the variables $\eta_0(X) + h_0(X)/\sqrt{n}$. If the distribution of $\eta_0(X)$ is smooth, then we might expect that the distribution functions of the perturbed variables $\eta_0(X) + h_0(X)/\sqrt{n}$ will be smooth as well. However, this appears not to be true in general, and it will usually be necessary to exploit some information about the functions $h_0$. (We need to consider functions in the support of the limit measure of the sequence $\sqrt{n}(\eta_n - \eta_0)$.) In this respect the conditions of Theorem 3.2 for composition with monotone functions are much more stringent than the conditions of Theorem 3.1 for the composition with Lipschitz functions.

The condition (19) is in terms of the indicator functions $1_{(-\infty, s]}$, and would have exactly the same form if we considered only indicator functions $\theta = 1_{(-\infty, \theta]}$, rather than general monotone functions. Thus the restrictive condition is connected to studying the classical empirical process.

The following lemma allows the verification of condition (19) in many cases. It will also be used in the next section to prove applicability of the delta-method. The lemma is similar to Lemma 5.1 of Ghoudi and Rémillard [6].

**Lemma 3.2.** *Suppose that $X, Y, Y_t$ (with $t > 0$) are real-valued random variables on a common probability space such that*

(i) *$X$ possesses a Lebesgue density $f$ that is continuous in a neighbourhood of $x$;*
(ii) *$\|Y_t - Y\|_\infty \to 0$ and $\|Y\|_\infty < \infty$;*
(iii) *the conditional distribution of $Y$ given $X = s$ can be represented by a Markov kernel $K(s, \cdot)$ such that the map $s \mapsto K(s, \cdot)$ is continuous at $x$ for the weak topology.*

*Then for every continuous function $g : \mathbb{R} \to \mathbb{R}$ and every converging sequences $x_t \to x$, $a_t \to a$ and $0 \le b_t \to b$, as $t \to 0$,*

$$\frac{1}{t} \mathrm{E} g(Y_t) 1_{x_t < X + t a_t Y_t \le x_t + t b_t} \to b \int g(y) \, K(x, dy) \, f(x).$$

*Proof.* First consider the case that $Y_t = Y$ for every $t$. By the definitions of $K$ and $f$, we can write

$$\frac{1}{t} \mathrm{E} g(Y) 1_{x_t < X + t a_t Y < x_t + t b_t}$$
$$= \frac{1}{t} \int \mathrm{E}\big(g(Y) 1_{(x_t - s)/t < a_t Y \le (x_t - s)/t + b_t} | X = s\big) f(s) \, ds$$
$$= \iint g(y) 1_{u < a_t y \le u + b_t} \, K(x_t - ut, dy) \, f(x_t - ut) \, du.$$

The inner integral is equal to $\mathrm{E} g(Y_t) 1_{u < a_t Y_t \le u + b_t}$ for $Y_t$ possessing the law $K(x_t - ut, \cdot)$. By assumption $a_t Y_t$ converges in distribution to the law of $aY$ for $Y$ possessing the law $K(x, \cdot)$. It follows that the inner integral converges to $\int g(y) 1_{u < ay < u + b} \, K(x, dy)$ for any $(u, b)$ such that $u$ and $u + b$ are not among the atoms of the law of $aY$. This includes almost every $u$ for every fixed $b$. Because $Y$ has bounded range, the double integral can be restricted $y$ in a compact set and hence $u$ in a compact set; the argument $x_t - ut$ fo $f$ is then restricted to a neighbourhood of $x$. Therefore, we can apply the dominated convergence theorem to see that the right side of the display converges to

$$\iint g(y) 1_{u < ay \le u + b} \, K(x, dy) \, f(x) \, du.$$

This reduces to $b \int g(y) \, K(x, dy) \, f(x)$ by Fubini's theorem. This concludes the proof of the lemma for the case that $Y_t = Y$.

Because $g$ is continuous, $Y$ possesses bounded range and $\|Y_t - Y\|_\infty \to 0$ we have that $\|g(Y_t) - g(Y)\|_\infty \to 0$. Therefore

$$\frac{1}{t} \mathrm{E} \big| g(Y_t) - g(Y) \big| 1_{x_t < X + t a_t Y_t < x_t + t b_t}$$
$$\le o(1) \frac{1}{t} \Pr\big(x_t - t a_t \|Y_t\|_\infty < X < x_t + t + t a_t \|Y_t\|_\infty\big).$$

This converges to zero as $X$ has a density that is bounded on bounded intervals.

Finally the difference $1_{x_t < X + t a_t Y_t < x_t + t b_t} - 1_{x_t < X + t a_t Y < x_t + t b_t}$ is nonzero only if $X + t a_t Y$ is in the union of intervals of total length bounded by $t a_t \|Y_t - Y\|_\infty$ in a neighbourhood of $x$. By the lemma with $Y_t = Y$ and $g = 1$, which is already proved, $t^{-1} \Pr(x_t < X + t a_t Y < x_t + t c_t) \to c f(x)$ for any sequences $x_t \to x$ and $c_t \to c$.

Hence this probability converges to zero for $c_t = ta_t\|Y_t - Y\|_\infty$, which satisfied $c_t \to 0$. We conclude from this that

$$\frac{1}{t}\big|\mathrm{E}g(Y)(1_{x_t < X + ta_t Y_t < x_t + tb_t} - 1_{x_t < X + ta_t Y < x_t + tb_t})\big|$$

$$\leq \|g(Y)\|_\infty \frac{1}{t}\mathrm{E}|1_{x_t < X + ta_t Y_t < x_t + tb_t} - 1_{x_t < X + ta_t Y < x_t + tb_t}|$$

converges to zero. The proof of the lemma is complete upon combining the preceding. □

In order to verify (19) with $K = \{h_0\}$ a single function we can apply the preceding lemma with $(X, Y)$ equal to the pair of variables $(\eta_0(X), h_0(X))$, $b_t = \delta_n$, and $t = 1/\sqrt{n}$. Then the conditions of the lemma require that the variable $\eta_0(X)$ possesses a continuous density, and that the conditional distribution of $h_0(X)$ given $\eta_0(X) = z$ depends continuously on the value of $z$. The second condition is clearly unpleasant, but appears to be natural in the present situation. It will involve a closer analysis of the support of the limit distribution of the sequence $\sqrt{n}(\eta_n - \eta_0)$.

To verify (19) with a general compact set $K \subset \ell^\infty(\mathcal{X}, \mathbb{R})$ we simply note that in view of the compactness it suffices to verify that for every sequence $h_n$ such that $\|h_n - h_0\|_\infty \to 0$ for some $h_0 \in H_0$,

$$\sup_{t \in \mathbb{R}^d, a \leq t \leq b} \sqrt{n}\, P\left(1_{\eta_0 + n^{-1/2}h_n \leq t + n^{-1/2}\delta_n} - 1_{\eta_0 + n^{-1/2}h_n \leq t}\right) \to 0.$$

Thus we can apply the preceding lemma with the variables $(X, Y, Y_t)$ equal to $(\eta_0(X), h_0(X), h_n(X))$ and $t = 1/\sqrt{n}$.

**Example 2, continued.** Suppose that $\eta_n$ is the empirical distribution of a random sample from the cumulative distribution function $\eta_0$ on $\mathbb{R}^d$. Then the limit distribution of the sequence $\sqrt{n}(\eta_n - \eta_0)$ is the $d$-dimensional $\eta_0$-Brownian sheet on $\mathbb{R}^d$.

If $d = 1$, then the Brownian sheet is a Brownian bridge and can be represented as $B \circ \eta_0$ for $B$ a standard Brownian bridge on the unit interval. A typical function in the support of the limit distribution of the sequence $\sqrt{n}(\eta_n - \eta_0)$ can be represented as $h_0 = h \circ \eta_0$ for some function $h : [0, 1] \to \mathbb{R}$. The conditional law of the variable $h_0(X)$ given $\eta_0(X) = z$ is the Dirac measure at $h(z)$. Because the standard Brownian bridge is continuous, the function $h$ can be taken continuous and hence the corresponding Markov kernels $K(z, \cdot) = \delta_{h(z)}(\cdot)$ are weakly continuous in $z$, as required by the preceding lemma.

If $d > 1$, then we can, without loss of generality, suppose that $\eta_0$ is a distribution function on $[0, 1]^d$ with uniform marginal distributions (i.e. a copula function). Then the conditioning event $\eta_0(X) = z$ will typically restrict $X$ to a one-dimensional curve in $[0, 1]^d$. Under sufficient smoothness of $\eta_0$, this curve will vary continuously with $z$, and under smoothness conditions on the law of $X$, the conditional distribution of $h_0(X)$ given $\eta_0(X) = z$ for a continuous function $h_0$ will vary continuously as well. Ghoudi and Remillard [6] give sufficient conditions for this continuity in a number of examples.

The preceding lemma can be extended to the case of multidimensional variables. For simplicity we only consider the two-dimensional case.

**Lemma 3.3.** *Suppose that $X, Y, Y_t$ (with $t > 0$) are random variables in $\mathbb{R}^2$ defined on a common probability space such that*

(i) $X$ possesses a Lebesgue density $f$ that has continuous conditional densities;

(ii) $\|Y_t - Y\|_\infty \to 0$ and $\|Y\|_\infty < \infty$;

(iii) the conditional distribution of $Y$ given $X = s$ can be represented by a Markov kernel $K(s, \cdot)$ such that the map $s \mapsto K(s, \cdot)$ is continuous at $x$ for the weak topology.

Then for every continuous function $g : \mathbb{R} \to \mathbb{R}$ and every converging sequences $x_t \to x$, $a_t \to a$ and $b_t \to b > 0$, as $t \to 0$,

$$\frac{1}{t} \mathrm{E} g(Y_t)(1_{X + t a_t Y_t \leq x_t + t b_t} - 1_{X + t a_t Y_t \leq x_t})$$

$$\to b_1 \int_{-\infty}^{x_2} \int g(y) \, K\big((x_1, s_2), dy\big) \, f(x_1, s_2) \, ds_2$$

$$+ b_2 \int_{-\infty}^{x_1} \int g(y) \, K\big((s_1, x_2), dy\big) \, f(s_1, x_2) \, ds_1.$$

*Proof.* The event $\{X + t a_t Y_t \leq x_t + t b_t\} \cap \{X + t a_t Y_t \leq x_t\}^c$ can be decomposed in the three events

$$I = \{x_{1t} < X_1 + t a_t Y_{1t} \leq x_{1t} + t b_{1t}, X_2 + t a_t Y_{2t} \leq x_{2t}\},$$

$$II = \{x_{1t} < X_1 + t a_t Y_{1t} \leq x_{1t} + t b_{1t}, x_{2t} < X_2 + t a_t Y_{2t} \leq x_{2t} + t b_{2t}\},$$

$$III = \{X_1 + t a_t Y_{1t} \leq x_{1t}, x_{2t} < X_2 + t a_t Y_{2t} \leq x_{2t} + t b_{2t}\}.$$

In view of the boundedness of the $Y_t$ the event $II$ is contained in an event of the form $\{X \in B_1\}$ for $B_t$ rectangles of area $O(t^2)$. Therefore, this event does not contribute to the limit.

The contribution of the event $I$ with $Y_t = Y$ can be written

$$\iiint g(y) 1_{u_1 < a_t y_1 \leq u_1 + b_{1t}} 1_{s_2 + t a_t y_{2t} \leq x_{2t}} K\big((x_{1t} - u_1 t, s_2), dy\big)$$

$$\times f(x_{1t} - u_1 t, s_2) \, du_1 \, ds_2.$$

By arguments as given previously this can be shown to converge to

$$\iiint g(y) 1_{u_1 < a y_1 \leq u_1 + b_1} K\big((x_1, s_2), dy\big) 1_{s_2 \leq x_2} f(x_1, s_2) \, du_1 \, ds_2.$$

The integral with respect to $u_1$ can be computed explicitly and the expression reduces to the first term on the right of the lemma.

The contribution of the event $III$ gives the second term.

We can replace $Y_t$ by $Y$ by similar arguments as in the one-dimensional case. (In fact, bound $x_{2t}$ by $\infty$ and use exactly the same arguments.) $\qquad\square$

## 4. Pseudo observations

In this section we consider the asymptotic behaviour of the process $\{\sqrt{n}(\mathbb{P}_n \theta \circ \eta_n - P\theta \circ \eta_0) : \theta \in \Theta\}$ for a given class $\Theta$ of functions $\theta : \mathbb{R}^d \to \mathbb{R}$. The set-up is the same as in Section 3. As explained in the introduction we can decompose this process as

$$\mathbb{G}_n(\theta \circ \eta_n - \theta \circ \eta_0) + \mathbb{G}_n \theta \circ \eta_0 + \sqrt{n} P(\theta \circ \eta_n - \theta \circ \eta_0).$$

Under the conditions of Theorem 3.1 or Theorem 3.2, Theorem 2.3, or their extensions, the first term will converge to zero in probability in $\ell^\infty(\Theta)$. The second term will converge in distribution to a Gaussian process in this space if and only if the class of functions $\Theta$ is Donsker for the law $P \circ \eta_0^{-1}$. If the third term also converges in distribution, then the sum of the three processes is asymptotically tight, and it will usually be straightforward to deduce its limit distribution from consideration of the marginal distributions.

The behaviour of the third term will follow by the (functional) delta-method if the sequence $\sqrt{n}(\eta_n - \eta_0)$ converges in distribution in the Banach space $H$ and the map $\eta \mapsto (P\theta \circ \eta : \theta \in \Theta)$ from $H$ to $\ell^\infty(\Theta)$ is suitably differentiable. If the limit distribution of the sequence $\sqrt{n}(\eta_n - \eta_0)$ concentrates on the space $H_0 \subset H$, then it suffices that the map $\eta \mapsto (P\theta \circ \eta : \theta \in \Theta)$ be "Hadamard differentiable tangentially to $H_0$", i.e. for every converging sequence $h_t \to h_0 \in H_0 \subset H$

$$\frac{1}{t} P\Big[ (\theta \circ (\eta_0 + t h_t) - \theta \circ \eta_0) \Big] \to L(h_0)(\theta),$$

uniformly in $\theta \in \Theta$, for a continuous linear map $L : \text{lin} \, H_0 \to \ell^\infty(\Theta)$. Under the additional condition that $L$ is defined on all $H$, this implies

$$\sqrt{n} P(\theta \circ \eta_n - \theta \circ \eta_0) = L(\sqrt{n}(\eta_n - \eta_0))(\theta) + o_P(1).$$

Cf. van der Vaart and Wellner [12], page 374.

As in the preceding section we consider the cases that the functions $\theta$ are smooth or of bounded variation separately. In the former case the differentiability is relative to a weak norm on $H$ (and is easy to prove), but for discontinuous functions $\theta$, such as the indicator functions $1_{(-\infty, \theta]}$, the differentiability requires a strong norm on $H$ and some conditions on the underlying distribution.

### 4.1. Smooth functions $\theta$

If the functions $\theta$ are differentiable with bounded derivatives, then the Hadamard differentiability is true for $H$ equipped with the $L_1(P)$-norm on $H$.

**Lemma 4.1.** *Let the functions $\theta : \mathbb{R}^d \to \mathbb{R}$ in $\Theta$ be continuously differentiable with derivative $\dot\theta$ such that $\|\dot\theta(x)\| \leq 1$ for every $x \in \mathbb{R}^d$. Then the map $\eta \mapsto (P\theta \circ \eta : \theta \in \Theta)$ from $L_1(\mathcal{X}, \mathcal{A}, P)$ to $\ell^\infty(\Theta)$ is Hadamard differentiable at $\eta_0$ with derivative $h \mapsto (P\dot\theta \circ h : \theta \in \Theta)$.*

*Proof.* Given a sequence $h_t$ with $P|h_t - h_0| \to 0$ we can write, by Fubini's theorem,

$$\Big| \frac{1}{t} P\big( \theta(\eta + t h_t) - \theta(\eta) \big) - P\dot\theta(\eta) h \Big|$$

$$= \Big| \int_0^1 P\big( \dot\theta(\eta + st h_t) h_t - \dot\theta(\eta) h \big) \, ds \Big|$$

$$\leq \int_0^1 P\big\| \dot\theta(\eta + st h_t) - \dot\theta(\eta) \big\| \|h_t\| \, ds + P\|\dot\theta(\eta)\| \|h_t - h_0\|.$$

The second term on the right is bounded above by $P|h_t - h_0|$ and converges to zero by assumption. The first term on the right converges to zero by the dominated convergence theorem. $\qquad \square$

### 4.2. Functions $\theta$ of bounded variation

In the second result we let $\Theta$ be a set of functions of bounded variation on a bounded interval in $\mathbb{R}$, and consider the Hadamard differentiability of the map $\eta \mapsto (P\theta \circ \eta : \theta \in \Theta)$ as a map from $\ell^\infty(\mathcal{X})$ to $\ell^\infty(\Theta)$. For simplicity of notation, let $X$ be a random variable with law $P$.

**Lemma 4.2.** *Let the functions $\theta \in \Theta$ be distribution functions of subprobability measures supported on a compact interval $I \subset \mathbb{R}$. Suppose that the variable $\eta_0(X)$ possesses a Lebesgue density $f$ that is continuous on a neighbourhood of $I$. Then the map $\eta \mapsto (P\theta \circ \eta : \theta \in \Theta)$ from $\ell^\infty(\mathcal{X})$ to $\ell^\infty(\Theta)$ is Hadamard differentiable at $\eta_0$ tangentially to the set of all $h_0$ such that there exists a version of the conditional distribution of $h_0(X)$ given $\eta_0(X) = s$ that is weakly continuous in $s \in I$. The derivative is given by $h_0 \mapsto \left( \int \mathrm{E}(h_0(X)|\eta_0(X) = s) f(s) \, d\theta(s) : \theta \in \Theta \right)$.*

*Proof.* Let $h_0$ be as given and suppose $h_t \to h_0$ in $\ell^\infty(\mathcal{X})$.

For given $s \in \mathbb{R}$ and $u > 0$ let $\chi_{s,u}$ be the continuous function that takes the value 0 on $(-\infty, s-u]$, takes the value 1 on $[s, \infty)$ and is linear on the interval $[s-u, s]$. Then $1_{[s,\infty)} \le \chi_{s,u}$, and hence

$$P(1_{s \le \eta_0 + th_t} - 1_{s \le \eta_0}) \le P(\chi_{s,u}(\eta_0 + th_t) - \chi_{s,u}(\eta_0)) + P(\chi_{s,u}(\eta_0) - 1_{[s,\infty)}(\eta_0)).$$

Because $\eta_0(X)$ possesses a Lebesgue density that is bounded on a neighbourhood of $I$ and $1_{[s,\infty)} - \chi_{s,u}$ vanishes off the set $(s-u, s)$, the second term on the right is bounded in absolute value by a multiple of $u$, uniformly in $s$ ranging through $I$, for small $u$. By choosing $u = \delta t$ this term divided by $t$ can be made arbitrarily small by choice of $\delta$.

Because $\chi_{s,u}$ is absolutely continuous with derivative $1/u$ on $(s-u, s)$ and 0 elsewhere, the first term on the right divided by $1/t$ can be written in the form

$$\int_0^1 \frac{1}{u} P(h_t 1_{s-u < \eta_0 + vth_t \le s}) \, dv.$$

For $u = \delta t$ this converges to $\mathrm{E}(h_0(X)|\eta_0(X) = s) f(s)$, by Lemma 3.2, uniformly in $s$ ranging over $I$.

It follows that, uniformly in $s$ ranging over $I$,

$$\limsup_{t \downarrow 0} \left( \frac{1}{t} P(1_{s \le \eta_0 + th_t} - 1_{s \le \eta_0}) - \mathrm{E}(h_0(X)|\eta_0(X) = s) f(s) \right) \le 0.$$

A similar argument using the functions $\chi_{s+u,u}$ instead of $\chi_{s,u}$ gives a corresponding lower bound, whence the expression in brackets converges to zero, uniformly in $s$ ranging through compacta. This concludes the proof of the lemma for $\Theta$ equal to the set of functions $1_{[s,\infty)}$ with $s$ in a compact interval.

For a general collection $\Theta$ of functions of bounded variation we can write

$$P(\theta(\eta_0 + th_t) - \theta(\eta_0)) = \int P(1_{s \le \eta_0 + th_t} - 1_{s \le \eta_0}) \, d\theta(s).$$

Next we use the assumption that the functions $\theta \in \Theta$ are supported on the compact interval $I$ with total variation bounded by 1.                                          $\square$

The applicability of the second lemma depends on whether the set $H_0$ of functions such that the conditional distribution of $h_0(X)$ given $\eta_0(X) = s$ is weakly continuous in $s$ is large enough to support the limit distribution of the sequence $\sqrt{n}(\eta_n - \eta_0)$. As noted in the preceding section, under some smoothness conditions on $\eta_0$ and on the distribution of $X$, the set $H_0$ typically contains all continuous functions. Then it suffices that the sequence possesses a continuous weak limiting process.

## 5. Examples: completion

In this section we return to two of the three examples discussed at the end of Section 3, Example 2 and Example 3. We give the theorems (and proofs) resulting from our approach. The general theme here is that the traditional results given in Corollaries 5.1 and 5.3 for indicator functions involve non-trivial restrictions on the underlying distribution $\eta_0$ of the data, while the results for indexing by Lipschitz functions given in Corollaries 5.2 and 5.4 involve almost no restrictions on $\eta_0$ (but significantly smoother indexing functions $\theta$).

### 5.1. Two corollaries for Kendall processes

For the Kendall process, Example 2, it suffices to consider the case in which $\eta_0$ is concentrated on $[0,1]^d$ and has uniform marginal distributions (i.e. is a copula function), as noted by Ghoudi and Remillard [6]. We first give a corollary for indexing by indicator functions, and then a corollary for indexing by Lipschitz functions.

**Corollary 5.1.** *Suppose that for a given interval* $[a,b] \subset (0,1)$:

 (i) *The variable* $\eta_0(X)$ *possesses a density* $k$ *with respect to Lebesgue measure that is continuous on a neighbourhood of* $[a,b]$.
 (ii) *The conditional distribution of* $X$ *given* $\eta_0(X) = s$, *has a regular version representable as a Markov kernel* $K(s,\cdot)$ *such that* $s \mapsto K(s,\cdot)$ *is continuous on* $[a,b]$ *for the weak topology.*

*Then the sequence of processes* $\sqrt{n}(K_n - K)$ *as in (11) tends in* $\ell^\infty[a,b]$ *in distribution to the process* $\left(\mathbb{G}_{\eta_0} f_\theta : \theta \in [a,b]\right)$ *for* $\mathbb{G}_{\eta_0}$ *an* $\eta_0$-*Brownian bridge process and* $f_\theta : [0,1]^d \to \mathbb{R}$ *defined as*

$$f_\theta(x) = 1_{\eta_0(x) \le \theta} - k(\theta)\, \mathrm{E}[1_{x \le X} | \eta_0(X) = \theta].$$

**Corollary 5.2.** *(Kendall processes, Example 2, indexed by Lipschitz functions). Suppose that* $\Theta$ *is a suitably measurable collection of continuously differentiable functions* $\theta : [0,1] \to [-1,1]$ *with derivatives* $\dot\theta$ *satisfying* $\|\dot\theta(x)\| \le 1$ *for every* $x \in [0,1]$. *Then the sequence of processes* $n^{-1/2} \sum_{i=1}^n \left(\theta(\eta_n(X_i)) - P\theta(\eta_0)\right)$ *tends in distribution in* $\ell^\infty(\Theta)$ *to the process* $\left(\mathbb{G}_{\eta_0} f_\theta : \theta \in \Theta\right)$ *for* $\mathbb{G}_{\eta_0}$ *an* $\eta_0$-*Brownian bridge process in* $\ell^\infty(\Theta)$ *and* $f_\theta : [0,1]^d \to \mathbb{R}$ *defined as*

$$f_\theta(x) = \theta\left(\eta_0(x)\right) - P\dot\theta(1_{x \le X}).$$

*Proof of Corollary 5.1.* We apply the decomposition (2) with $f_{\theta,\eta}x = 1\{\eta(x) \le \theta\}$, for distribution functions $\eta$ on $[0,1]^d$, $\theta \in [0,1]$ and $x \in [0,1]^d$.

As discussed following the proof of Lemma 3.2, hypotheses (i) and (ii) imply that the condition (19) for Theorem 3.2 (with $d = 1$) holds by way of Lemma 3.2, and hence the first term on the right side of (2) converges in probability to 0 uniformly in $\theta \in [a,b]$.

The second term is simply the usual empirical process for the i.i.d. one-dimensional random variables $\eta_0(X_1), \ldots, \eta_0(X_n)$, and hence it converges weakly as claimed by standard theory.

To handle the third term, note that (i) and (ii) imply that the hypotheses of Lemma 4.2 hold, and hence that the map $\eta \mapsto \{Pf_{\theta,\eta} : \ \theta \in \Theta\}$ from $\ell^\infty(\mathcal{X})$ to $\ell^\infty([a,b])$ is Hadamard differentiable tangentially to $C([0,1]^d)$ with derivative $L : C([0,1]^d) \to \ell^\infty([a,b])$ given by

$$L(h_0)(\theta) = -\mathrm{E}\big(h_0(X)\|\eta_0(X) = \theta\big)k(\theta).$$

Weak convergence of the third term then follows from van der Vaart and Wellner [12], Theorem 3.9.5, page 375.

The joint limit law of the second and third term can be determined from the marginals, and the limit of the sum of the two terms can be represented in the form as given. An insightful way to derive this is from asymptotic linearity of the two terms as follows. The second term is already linear with influence functions $x \mapsto 1_{\eta_0(x) \le \theta}$. The third term can be approximated by $L((\sqrt{n}(\eta_n - \eta))$, where $\eta_n - \eta = n^{-1}\sum_{i=1}^{n}(1_{X_i \le x} - \eta(x))$, so that $L((\sqrt{n}(\eta_n - \eta)) = n^{-1/2}\sum_{i=1}^{n} L(1_{[X_i,1]} - \eta)$. The terms in the latter sum should be understood as $L$ acting on the functions $x \mapsto 1_{[X_i,1]}(x) - \eta(x)$ for fixed $X_i$. We thus obtain that

$$L((\sqrt{n}(\eta_n - h)) = n^{-1/2}\sum_{i=1}^{n} L(1_{[X_i,1]}) - \sqrt{n}L(\eta) = \mathbb{G}_{\eta_0} L(1_{[X_i,1]}).$$

The representation of the limit process as given in the corollary follows. $\qquad\blacksquare$

For many distribution functions $\eta_0$ the corresponding density $k$ of $K$ is unbounded at 0 and hence not continuous on $[0,1]$. See Barbe, et al. [2], pages 202-208, for a number of explicit examples. In particular this is true even when $\eta_0$ is the uniform distribution on $[0,1]^d$. For such distributions the preceding corollary does not yield convergence of Kendall's process in the space $\ell^\infty([0,1])$. However, this convergence may be valid even when $k$ is unbounded. Barbe et al. [2] show that under the growth condition

$$k(t) = o(t^{-1/2}(\log(1/t))^{-1/2-\epsilon}), \qquad t \downarrow 0, \quad \epsilon > 0.$$

convergence in the full domain still holds. They achieve this using results of Alexander [1] to show that the empirical process $\sqrt{n}(\eta_n - \eta_0)1\{\eta_0 \ge a_n\}$ converges in the weighted metric $\|\cdot/q(\eta_0)\|_\infty$ if $q(t) = t^{1/2}(\log(1/t))^p$ for some $1/2 < p < r/2$ and $a_n = n^{-1}(\log n)^r$. This strengthening of the convergence of $\sqrt{n}(\eta_n - \eta_0)$ then compensates for the growth of $k$ at 0.

*Proof of Corollary 5.2.* This follows by combining Theorems 3.1 and Lemma 4.1 with the fact that $\mathcal{F} = \{\theta \circ \eta_0 : \theta \in \Theta\}$ is Donsker. $\qquad\blacksquare$

### 5.2. *Two corollaries for copula processes*

For the copula processes (12) in Example 3 it again suffices to consider the case in which $\eta_0 = C$, so that all all the marginal distributions $\eta_{0,j}$, $j = 1, \ldots, d$, are Uniform$(0,1)$. The first of the following two corollaries was obtained in Stute [9] and Ghoudi and Remillard [7].

**Corollary 5.3.** *Suppose that:*

  (i) *$\eta_0 = C$ is continuous.*

  (ii) *The copula function $\eta_0 = C$ is continuously differentiable on $[0,1]^d$ with gradient $\nabla C(u)$.*

*Then the sequence of copula processes $\sqrt{n}(C_n - C)$ given in (12) converges in distribution in $\ell^\infty([0,1]^d)$ to the process $(\mathbb{G}_{\eta_0} f_u : u \in [0,1]^d)$ for $\mathbb{G}_{\eta_0}$ an $\eta_0$-Brownian bridge process, and $f_u : [0,1]^d \to \mathbb{R}$ defined as*

$$f_u(x) = 1_{x \le u} - \nabla C(u)'(1_{x_1 \le u_1}, \dots, 1_{x_d \le u_d}).$$

**Corollary 5.4.** *(Copula processes, Example 3, indexed by Lipschitz functions). Suppose that $\Theta$ is a suitably measurable collection of continuously differentiable functions $\theta : [0,1]^d \to \mathbb{R}$ such that with derivative $\|\dot{\theta}(x)\| \le 1$ for every $x \in [0,1]^d$ and satisfying $J(1, \Theta, L_2) < \infty$. Then the sequence of processes $n^{-1/2} \sum_{i=1}^n (\theta(\eta_n(X_i)) - P\theta(\eta_0))$ tends in distribution in $\ell^\infty(\Theta)$ to the process $(\mathbb{G}_{\eta_0} f_\theta : \theta \in \Theta)$ for $\mathbb{G}_{\eta_0}$ an $\eta_0$-Brownian bridge process in $\ell^\infty(\Theta)$ and $f_\theta : [0,1]^d \to \mathbb{R}$ defined as*

$$f_\theta(x) = \theta(x) - P\dot{\theta}(1_{x_1 \le X_1}, \dots, 1_{x_d \le X_d}).$$

*Proof of Corollary 5.3.* We apply the decomposition (2) with

$$f_{\theta,\eta}(x) = 1\{\eta_1(x_1) \le \theta_1, \dots, \eta_d(x_d) \le \theta_d\}$$

for $\theta = (\theta_1, \dots, \theta_d) \in [0,1]^d$, $x = (x_1, \dots, x_d) \in [0,1]^d$, and $\eta_j$ the $j$th one-dimensional marginal distribution function on $[0,1]$ of the distribution function $\eta$ (so $\eta_j(u_j) = \eta(1, \dots, 1, u_j, 1, \dots, 1)$).

To show that the first term in (2) converges to zero uniformly in $\theta \in [0,1]^d$, we can apply Theorem 2.1. The class of functions

$$f_{\theta,\eta}(x) = 1\{x_1 \le \eta_1^{-1}(\theta_1), \dots, x_d \le \eta_d^{-1}(\theta_d)\}$$

is a class of indicators of a Vapnik-Chervonenkis-class of sets. Thus Theorem 2.1 applies if we show that (3) holds. But this is easily verified by the assumed continuity of $\eta_0 = C$ and the uniform consistency of the empirical quantile functions $\eta_{n,j}^{-1}$ for $j = 1, \dots, d$. Thus (1) holds.

The second term in (2) is simply the classical empirical process of the random vectors $X_1, \dots, X_n$ in $[0,1]^d$, and converges weakly by classical theory.

Finally, the third term in (2) converges weakly to $\nabla C(u)' \cdot \mathbb{G}_{\eta_0}(v(X, u))$, for $v(x, u) = (1_{x_1 \le u_1}, \dots, 1_{x_d \le u_d})$, by the delta-method for the map $\eta \mapsto P f_{u,\eta}$. This map can be decomposed as

$$\eta \mapsto (\eta_1^{-1}(u_1), \dots, \eta_d^{-1}(u_d)) \mapsto (C \circ (\eta_1^{-1}(u_1), \dots, \eta_d^{-1}(u_d)), \qquad u \in [0,1]^d$$

and can be shown to Hadamard-differentiable from the domain of distribution functions in $\ell^\infty(\mathcal{X}) = \ell^\infty([0,1]^d)$ to $\ell^\infty(\Theta) = \ell^\infty([0,1]^d)$ by the chain rule, using the continuity of $\nabla C$ and the fact that the quantile transformation is Hadamard differentiable. $\square$

It is possible to extend Corollary 5.3 to the case in which $\nabla C$ is continuous on $(0,1)^d$ but satisfies certain growth restrictions at 0 and/or 1. Then weighted metrics are involved in the proof.

*Proof of Corollary 5.4.* This follows by combining Theorems 3.1 and Lemma 4.1 with the fact that $\mathcal{F} = \{\theta : \theta \in \Theta\}$ is Donsker and the delta-method, e.g. van der Vaart and Wellner [12], Theorem 3.9.5, page 375. $\square$